\documentclass{amsart}

\usepackage{amssymb,latexsym}

\newcommand{\doublesubscript}[3]{
\displaystyle\mathop{\displaystyle #1_{#2}}_{#3} }

\numberwithin{equation}{section}
\newtheorem{theorem}{Theorem}[section]
\newtheorem{proposition}[theorem]{Proposition}
\newtheorem{corollary}[theorem]{Corollary}

\newtheorem{remark}[theorem]{Remark}
\newtheorem{lemma}[theorem]{Lemma}

\newtheorem{example}[theorem]{Example}

\newcommand{\mat}[4]{\left(\!\!\begin{array}{cc}
#1 & #2 \\ #3 & #4 \\
\end{array}\!\!\right)}

\def\ZZ{\mathbb{Z}}
\def\CC{\mathbb{C}}
\def\gg{\mathfrak{g}}
\def\hh{\mathfrak{h}}
\def\nn{\mathfrak{n}}
\def\ii{\mathbf{i}}

\begin{document}

\title
{On symplectic leaves and integrable systems in standard
complex semisimple Poisson-Lie groups}

\author {Mikhail Kogan}

\author{Andrei Zelevinsky}

\begin{abstract}
We provide an explicit description of symplectic leaves of a
simply connected connected semisimple complex Lie group equipped
with the standard Poisson-Lie structure. This sharpens previously
known descriptions of the symplectic leaves as connected components
of certain varieties. Our main tool is the machinery of twisted generalized minors.
They also allow us to present several quasi-commuting coordinate systems
on every symplectic leaf.
As a consequence, we construct new completely integrable systems on
some special symplectic leaves.
\end{abstract}

\address{\noindent Department of Mathematics, Northeastern University,
  Boston, MA 02115}
\email{misha@neu.edu}

\email{andrei@neu.edu}

\date{March 7, 2002}

\thanks{Research supported in part by NSF grants:
Postdoctoral Research Fellowship (M.K.)
and DMS-9971362 (A.Z.).}

\maketitle

\pagestyle{myheadings} \markboth{MIKHAIL KOGAN AND ANDREI
ZELEVINSKY}{SYMPLECTIC LEAVES AND INTEGRABLE SYSTEMS}

\section{Introduction}

Let $G$ be a simply connected connected semisimple complex Lie
group supplied with the standard Poisson-Lie structure. It is well
known (see e.g., \cite{CC,hl,r+}) that the symplectic leaves in
$G$ are closely related to \emph{double Bruhat cells}, the
intersections of double cosets of two opposite Borel subgroups in
$G$. Double Bruhat cells were studied in \cite{fz}; one of the
main tools developed there was a family of regular functions on
them called \emph{twisted (generalized) minors}. In the present
paper, we provide some applications of these functions to the
study of symplectic leaves in $G$ and integrable systems on them
(some results in this direction were obtained in \cite{r+,r,y}).
One of our main goals is to bring the machinery of twisted minors
to the attention of the experts in the field. We believe these
functions should have further applications to the study of
symplectic leaves and integrable systems.

It was shown in \cite{fz} that twisted minors give rise to a
family of toric charts in every double Bruhat cell. These charts
were used in \cite{z-imrn} for determining the connected
components of real double Bruhat cells. In a similar spirit, we
obtain here an explicit description of symplectic leaves in $G$
(Theorem~\ref{symplectic-leaves}); furthermore, using twisted
minors we produce toric charts in every symplectic leaf. This
sharpens the results in \cite{hl,r+,y}, where the symplectic
leaves were characterized as connected components of some
subvarieties of the double Bruhat~cells.

Our second main result (Theorem~\ref{poisson-brackets}) asserts
that certain twisted minors quasi-commute with each other (that
is, their Poisson bracket is a scalar multiple of their product),
and gives an explicit expression for their Poisson bracket. (In a
different context, this calculation will also appear in a
forthcoming paper by A.~Berenstein and one of the authors (A.Z.)
devoted to the study of \emph{quantum} double cells.)

As an application of Theorem~\ref{poisson-brackets}, we show in
Corollary~\ref{cor:integrable-system} that the twisted minors can
be used to construct integrable systems on some special symplectic
leaves. We are unable to match these integrable systems with any
known ones; we call for experts to try to recognize them.
It should be possible to extend the method of
Corollary~\ref{cor:integrable-system} to produce more examples of
integrable systems; we believe that this method deserves further study.

The paper is organized as follows.
In Section~\ref{main-results}, we provide necessary background
and state our main results.
Their proofs are given in Sections~\ref{sec:proofs-symplectic}
and \ref{sec:proofs-integrable}.

\smallskip

\textsc{Acknowledgments.} We thank Kolya Reshetikhin for patiently explaining to one of us
(A.Z.) the basics of symplectic leaves, and Arkady Berenstein whose
comments and suggestions helped us to simplify our original proof
of Theorem~\ref{poisson-brackets}. We also thank Misha Gekhtman
and Misha Shapiro for helpful discussions and for informing us
about some of their yet unpublished results obtained jointly with Alek Vainshtein.

\section{Main results}
\label{main-results}

\subsection{Poisson manifolds and symplectic leaves}
We start by recalling some basic definitions.
Let $M$ be a smooth manifold.
Denote by $C^\infty (M)$ the set of smooth complex
valued functions on $M$.
Recall that $M$ is a \emph{Poisson manifold} if it is equipped
with a bilinear map $\{\cdot ,\cdot\}:C^\infty (M) \times
C^\infty(M)\to C^\infty (M)$ called {\it Poisson bracket}, which
makes $C^\infty(M)$ into a Lie algebra and satisfies the Leibniz
identity
$\{fg,h\}=f\{g,h\}+\{f,h\}g$
for $f,g,h\in C^\infty(M)$.
In all the cases we consider $M$ will be a complex algebraic
variety; then a Poisson bracket extends uniquely to a
Poisson bracket on the space of rational functions on $M$.

A smooth map $\varphi:M\to N$ between two Poisson manifolds is called
\emph{Poisson} if $\varphi^*\{f,g\}_N=\{\varphi^*(f),\varphi^*(g)\}_M$ for
every $f,g\in C^\infty(N)$;
here $\{.,.\}_M$ and $\{.,.\}_N$
are the Poisson brackets on $M$ and $N$, respectively.
The Poisson structure on the product
$M\times N$ of two Poisson manifolds is defined by
$$
\{f,g\}(x,y) = \{f(.,y),g(.,y)\}_M (x) +  \{f(x,.),g(x,.)\}_N (y)
$$
for $f,g\in C^\infty(M\times N)$.

For a function $f$ on a Poisson manifold $M$, the
\emph{Hamiltonian vector field} $X_f$ is the vector field
associated to the derivation $\{f,\cdot\}$. For a point $p\in M$,
the restrictions of Hamiltonian vector fields to $p$ form a vector
subspace $V_p$ inside the tangent vector space $T_pM$. This
defines a distribution $V$ on $M$. This distribution is known to
be integrable (see e.g. \cite{ks}). Thus, every Poisson manifold
$M$ is a disjoint union of connected symplectic manifolds
$S_\alpha$, such that $T_p S_\alpha=V_p$ for every $p\in S_\alpha$
and the symplectic form $\omega_\alpha$ on $S_\alpha$ is given by
$$
\omega_\alpha(X_f,X_g)(p)=\{f,g\}(p){\ \rm { for }}\ p\in S_\alpha, \
f,g\in C^\infty(M).
$$
The manifolds $S_\alpha$ are called the \emph{symplectic leaves}
of $M$. Clearly, every inclusion $S_\alpha\to M$ is Poisson.

\subsection{Standard Poisson-Lie structure on a semisimple Lie
group}

Recall that a \emph{Poisson-Lie group} is a Lie group $G$ equipped
with a Poisson bracket
such that the multiplication map $G\times G\to G$ is Poisson.

\begin{example}
\label{example} {\rm Let
$G=SL_2=\{\mat{x_{11}}{x_{12}}{x_{21}}{x_{22}}\ : \
x_{11}x_{22}-x_{12}x_{21}=1\}$; then we can define (see e.g.
\cite{ks}) a family of Poisson structures on $G$ parameterized by
$d\in \CC$~by
$$
\begin{array}{ccc}
\{x_{12},x_{11}\}= d x_{11} x_{12}, & \{x_{21},x_{11}\}= d x_{11}
x_{21}, &
\{x_{22},x_{11}\}= 2d x_{12} x_{21}, \\
\{x_{12},x_{21}\}=0, & \{x_{22},x_{12}\}= d x_{12} x_{22}, &
\{x_{22},x_{21}\}= d x_{21}x_{22} \ .
\end{array}
$$
We will indicate that $SL_2$ is equipped with the above Poisson
structure by writing~$SL_2^{(d)}$.}
\end{example}

Let $G$ be a simply-connected connected semisimple complex Lie
group with the Lie algebra~$\gg$. Let $f_i, \alpha_i^\vee, e_i \,
(i = 1, \dots, r)$ be the Chevalley generators of $\gg$, and
$\gg=\nn_- \oplus \hh \oplus \nn$ be the corresponding triangular
decomposition. Let $\alpha_1,\dots,\alpha_r \in \hh^*$ be the
simple roots of $\gg$. The Cartan matrix $A=(a_{ij})$ is given by
$a_{ij}= \langle \alpha_i^\vee, \alpha_j\rangle$. We fix a
diagonal matrix $D$ with positive diagonal entries
$d_1,\dots,d_r$, which symmetrizes~$A$, i.e., $d_i a_{ij} = d_j
a_{ji}$ for all $i$ and $j$. For $t\in\CC$, define
$$
x_i(t)=\exp(te_i), \ \ x_{\bar i} (t)=\exp(tf_i) \ .
$$
The canonical inclusions $\varphi_i: SL_2\to G$ are defined by
$$
\varphi_i \mat{1}{t}{0}{1}=x_i(t), \  \ \varphi_i
\mat{1}{0}{t}{1}=x_{\bar i}(t).
$$

The \emph{standard Poisson-Lie structure} on $G$ is uniquely
determined by the requirement that every map $\varphi_i:
SL_2^{(d_i)}\to G$ is Poisson. (The uniqueness is easy to see;
an explicit construction using Manin triples can be found in \cite{ks} or \cite{r}.)

\subsection{Generalized minors and their twists}

To state our main results
we need to introduce more notation and recall some results of \cite{fz}.
Let $N_-, H, N$ be the subgroups of $G$, which correspond to
$\nn_-, \hh$ and $\nn$. We set $B_-=HN_-$ and $B=HN$ to be the
pair of opposite Borel subgroups. The set $G_0=N_-HN$ is a Zariski
open subset of $G$ consisting of the elements $x\in G$, which have
(a unique) Gaussian decomposition $x=[x]_-[x]_0 [x]_+$, with
$[x]_- \in N_-, [x]_0 \in H$, and $[x]_+ \in N$.

The \emph{weight lattice} $P$ consists of elements
$\gamma\in\hh^*$ such that $\langle \alpha_i^\vee, \gamma \rangle
\in \ZZ$ for all~$i$. Every weight $\gamma\in P$ defines a
multiplicative character $a\mapsto a^\gamma$ of $H$, defined by
$\exp(h)^\gamma= e^{\langle h, \gamma \rangle}$ for $h\in\hh$. The
basis of \emph{fundamental weights} $\omega_1,...,\omega_r$ in $P$
is defined by $\langle \alpha_i^\vee, \omega_j \rangle =
\delta_{ij}$.

The Weyl group $W$ is defined by $W = {\text {Norm}}_G(H)/H$.
It acts on $H$ by $a^w=w^{-1}aw$, and the following formula
defines the action of $W$ on $P$:
$$
a^{w(\gamma)}=(w^{-1} a w)^\gamma,{\ \rm { for }}\ a\in H, w\in W,
\gamma\in P.
$$
The Weyl group is generated by simple reflections $s_1,...,s_r$
acting on weights by $s_i(\gamma)=\gamma- \langle \alpha_i^\vee, \gamma
\rangle \alpha_i$. If $w=s_{i_1}...s_{i_m}$ is a shortest possible
expression of $w$ as a product of simple reflections, then
$(i_1,...,i_m)$ is called a \emph{reduced word} of $w$, and $m$ is
called the length of $w$ and denoted by $\ell(w)$.

For every $w\in W$, we define a special representative $\overline
w$ of $w$ in ${\text {Norm}}_G(H)$ as follows. Let
$$
\overline s_i=\varphi_i\mat{0}{-1}{1}{0}.
$$
Then, if we also require $\overline w_1\overline w_2
=\overline{w_1w_2}$ as long as $\ell(w_1)+\ell(w_2)=\ell(w_1w_2)$,
it is not difficult to see that $\overline w$ is well defined for
every $w\in W$ (see \cite{fz}).

For $u,v\in W$, define {\it the double Bruhat cell $G^{u,v}$} to
be
$$G^{u,v} = BuB\cap B_-vB_- \ .$$
By \cite[Theorem~1.1]{fz}, $G^{u,v}$ is a smooth algebraic variety
of dimension $\ell (u) + \ell (v) + r$ (recall that $r$ is the
rank of $\gg$).

For $x\in G_0$ and a fundamental weight $\omega_i$, define
$$
\Delta_i(x)=[x]_0^{\omega_i}.
$$
It is shown in \cite{fz} that $\Delta_i$ extends to a regular
function on $G$. For type $A_r$ (when $G = SL_{r+1}$), this is
just the principal $i \times i$ minor of a matrix $x$.

For any pair $u,v\in W$, the corresponding \emph{generalized
minor} is a regular function on $G$ given by
$$
\Delta_{u\omega_i,v\omega_i} (x) =\Delta_i({\overline u}^{\ -1}x
\overline v).
$$
It is shown in \cite{fz} that these functions are well defined,
that is they depend only on the weights $u\omega_i$ and
$v\omega_i$ and do not depend on the particular choice of $u$ and
$v$.

Define an involutive automorphism $x\mapsto x^\theta$ on $G$ by
$$
a^\theta=a^{-1} \ \ (a\in H),\ \  x_i(t)^\theta=x_{\bar i}(t),\  \
x_{\bar i}(t)^\theta=x_i(t).
$$
The \emph{twist map} for $u,v\in W$ is a biregular isomorphism
$x\mapsto x'$ between $G^{u,v}$ and $G^{u^{-1},v^{-1}}$ given by
(see \cite[Theorem~1.6]{fz})
\begin{equation}
\label{eq:twist} x'=([{\overline u}^{\ -1}x]_-^{-1}{\overline
u}^{\ -1}x \overline{v^{-1}} [x\overline{v^{-1}}]_+^{-1})^\theta.
\end{equation}

A \emph{double reduced word} of $(u,v)$ is a word ${\bf
i}=(i_1,\dots,i_m)$ of length $m = \ell (u) + \ell (v)$ in the
alphabet $[1,\dots,r] \cup [\bar 1,\dots,\bar r]$ such that the
subword of ${\bf i}$ consisting of all letters from $[\bar
1,\dots,\bar r]$ is a reduced word of $u$, and the subword
consisting of all letters from $[1,\dots,r]$ is a reduced word of
$v$. \
For $i=1,\dots,r$, we denote $\varepsilon(i)=+1$ and
$\varepsilon(\bar i)=-1$, and set $|i|=|\bar i|=i$.

In what follows, we fix $u,v \in W$ and a double reduced word
${\bf i}$ of $(u,v)$. We append $r$ entries
$i_{m+1},\dots,i_{m+r}$ to $\bf i$ by setting $i_{m+j}=\bar j$.
For $k = 1, \dots, m$, we set
$$
u_{\geq k}=\doublesubscript{\prod}{\ell=m,\dots,k}
{\varepsilon(i_\ell) = -1} s_{|i_\ell|}, \ \
v_{<k}=\doublesubscript{\prod}{\ell=1,\dots,k-1}
{\varepsilon(i_\ell) = +1} s_{|i_\ell|},
$$
where the notation implies that the index $\ell$ in the first
(resp.~second) product is decreasing (resp.~increasing). We also
set $u_{\geq k}=e$, $v_{<k}=v$ for $k = m+1, \dots, m+r$. For
example, if ${\bf i}=(1,\bar 2, 2, \bar 3, 3, 2,\bar 1)$ then
$u_{\geq 4}=s_1s_3$, $v_{<4}=s_1s_2$.

For every $k = 1, \dots, m+r$, we set
$$\gamma^k= u_{\geq k}\omega_{|i_k|}, \, \, \delta^k = v_{<k}\omega_{|i_k|}$$
and introduce a regular function $M_k$ on $G^{u,v}$ by setting
\begin{equation}
\label{eq:twisted-minors-def} M_k (x) =\Delta_{\gamma^k,
\delta^k}(x'),
\end{equation}
where $x'$ is given by (\ref{eq:twist}).
We refer to the family $M_1, \dots, M_{m+r}$ as \emph{twisted
minors} associated with a reduced word $\bf i$. Their significance
stems from the following result (see \cite[Theorems~1.2, 1.9, 1.10
and formula (1.21)]{fz}).

\begin{theorem}
\label{biregular} The map $x_{\bf i}:H\times \CC^m \to G$ given by
$$
x_{\bf i}(a;t_1,\dots,t_m)=ax_{i_1}(t_1)\dots x_{i_m}(t_m)
$$
restricts to a biregular isomorphism between a complex torus
$H\times (\CC - \{0\})^m$ and a Zariski open subset $U_\ii=\{x \in
G^{u,v} : M_k (x) \neq 0 \text { for } 1 \leq k \leq m+r\}$ of the
double Bruhat cell $G^{u,v}$. Furthermore, for $k = 1, \dots, m+r$
and $x = x_{\ii}(a;t_1,\dots,t_m) \in U_\ii$, we have
\begin{equation}
\label{eq:M-k-monomial} M_k(x)=a^{-u \gamma^k}
\doublesubscript{\prod}{1 \leq \ell<k}{\varepsilon(i_\ell) = -1}
t_{\ell}^{\langle \alpha^\vee_{|i_\ell|},
u^{-1}_{\geq\ell}\gamma^k\rangle} \doublesubscript{\prod}{k \leq
\ell \leq m}{\varepsilon(i_\ell) = +1} t_{\ell}^{\langle
\alpha^\vee_{|i_\ell|},v^{-1}_{<\ell+1}\delta^k \rangle}.
\end{equation}

\end{theorem}

\subsection{Symplectic leaves of the standard  Poisson-Lie structure}

We are now ready to  describe symplectic leaves of the standard
Poisson-Lie structure on $G$. We retain the above terminology and
notation.

Given $u,v\in W$, let $H^{u,v}$ denote the subtorus of $H$ formed
by elements $(a^u)^{-1}\cdot a^v$ for $a\in H$. We also  denote by
$I(u,v)$ the set of all indices $i$ such that
$u\omega_i=v\omega_i=\omega_i$.

\begin{theorem}
\label{symplectic-leaves}
For every $u,v\in W$, the set
$$S^{u,v}=\{x\in G^{u,v}  :  [{\overline u}^{\ -1}x]_0 \cdot ([x
\overline {v^{-1}}]_0)^v \in H^{u,v},\ [{\overline u}^{\ -1}x]_0^{\omega_i}=1
\ {\rm \ for \ all \ } i \in I(u,v)\}$$
is a symplectic leaf in $G$.  Every symplectic leaf in $G$ is of the form
$S^{u,v} \cdot a$ for some $u,v\in W$ and $a \in H$.
\end{theorem}

In the process of the proof of Theorem~\ref{symplectic-leaves}, we
show that every double reduced word of $(u,v)$
gives rise to a toric chart in $S^{u,v}$.
More precisely, we prove the following.

\begin{proposition}
\label{pr:Suv-toric-chart}
For every double reduced word
$\ii =(i_1,\dots,i_m)$ of $(u,v)$, the intersection
$S^{u,v}\cap U_\ii$ is a dense Zariski open
subset of $S^{u,v}$, and the mapping
$$x\mapsto(M_1(x),\dots,
M_m(x),[{\overline u}^{\ -1}x]_0 \cdot ([x \overline
{v^{-1}}]_0)^v)$$
is a biregular isomorphism between
$S^{u,v}\cap U_\ii$ and the complex torus $(\CC-\{0\})^m\times H^{u,v}$.
In particular, the symplectic leaf $S^{u,v}$ has (complex) dimension
$m + \dim H^{u,v}$, where $m = \ell (u) + \ell (v)$.
\end{proposition}

Note that the set
\begin{equation}
\label{eq:tildeS}
\tilde S^{u,v}=\{x\in G^{u,v} : [{\overline
u}^{\ -1}x]_0 \cdot ([x \overline {v^{-1}}]_0)^v
 \in H^{u,v}\}
\end{equation}
has been known to be a union of finite number of symplectic leaves
(see \cite{hl,r,y} and
Proposition~\ref{pr:symplectic-leaves-thru-tildeS} below). The new
result in Theorem~\ref{symplectic-leaves} is an explicit
description of the connected components of $\tilde S^{u,v}$.
As a consequence of Theorem~\ref{symplectic-leaves}, we obtain the
following corollary.

\begin{corollary}
\label{cor:number-of-components}
The number of connected components of $\tilde S^{u,v}$ is equal to $2^{|I(u,v)|}$.
\end{corollary}

\subsection{Poisson brackets of twisted minors and integrable systems on
special symplectic leaves}

Let $(\gamma, \gamma')$ denote the
$W$-invariant scalar product on $\hh^*$ such that $(\alpha_i,
\gamma)=d_i \langle \alpha_i^\vee, \gamma \rangle$ for all $i$ and~$\gamma$.
Assume we are given $u,v\in W$ and a double reduced
word $\ii$ of $(u,v)$, and let $M_1,\dots M_{m+r}$ be
twisted minors associated with $\ii$.
Here is our next main result.

\begin{theorem}
\label{poisson-brackets} On every symplectic leaf in a double
Bruhat cell $G^{u,v}$, the standard Poisson bracket between
twisted minors is given by
\begin{equation}
\label{eq:M-poisson} \{ M_k,M_{k'}\}= ((\gamma^k,\gamma^{k'})-
(\delta^k,\delta^{k'}))M_k M_{k'}
\end{equation}
for $1 \leq k \leq  k' \leq m+r$.
\end{theorem}

Let $u$ be an arbitrary element of the Weyl group $W$. Using
Theorem~\ref{poisson-brackets}, we shall construct a family of
completely integrable systems on the symplectic leaf $S^{u,u}$,
one for each reduced word ${\bf j} = (j_1,\dots,j_{\ell(u)})$ of
$u$. To do this, we first associate with~$\bf j$ a double reduced
word ${\bf i} = (i_1, \dots, i_m)$ of $(u,u)$ (where $m = 2 \ell (u)$)
by $i_{2k-1}=j_k$ and $i_{2k}=\overline {j_k}$ for $1\leq
k\leq \ell(u)$. Let $M_1, \dots, M_{m+r}$ be twisted minors
associated with $\bf i$.

\begin{corollary}
\label{cor:integrable-system} The twisted minors $M_{2k-1}$ for $k
= 1, \dots, \ell (u)$ form a completely integrable system
on the symplectic leaf $S^{u,u}$, i.e. they are independent on
$S^{u,u}$, Poisson commute with each other, and the cardinality
$\ell (u)$ of this family is equal to $\frac{1}{2} \dim S^{u,u}$.
\end{corollary}

\begin{remark}
\label{remark} {\rm M.~Gekhtman and M.~Shapiro informed us that
they have proved (in an ongoing joint work with A.~Vainshtein)
that the twist isomorphism (\ref{eq:twist}) between double cells $G^{u,v}$ and
$G^{u^{-1},v^{-1}}$ is an anti-isomorphism of Poisson manifolds
(this means that it becomes Poisson if we change the sign of the
Poisson bracket on $G^{u^{-1},v^{-1}}$). Thus, Corollary
\ref{cor:integrable-system} also gives rise to a family of integrable
systems on symplectic leaves $S^{u,u}$ formed by (non-twisted)
generalized minors.}
\end{remark}

\begin{example}{\rm
Let $u = v = w_0$ be the longest element in $W$. The corresponding
double Bruhat cell $G^{w_0, w_0}$ is an open set in $G$ given by
$$G^{w_0, w_0} = \{x \in G :
\Delta_{w_0 \omega_i, \omega_i}(x) \neq 0, \Delta_{\omega_i,w_0
\omega_i}(x) \neq 0 {\rm \ for \ all \ } i\} \ .$$
Let $i \mapsto i^*$ be an involution on the index set $\{1, \dots, r\}$
induced by the action of $(- w_0)$ on fundamental weights: that is, we have
$w_0 (\omega_i) = - \omega_{i^*}$.
As an easy consequence of Theorem~\ref{symplectic-leaves},
the symplectic leaf $S^{w_0, w_0}$ is given by
$$S^{w_0, w_0} = \{x \in G^{w_0,w_0}:
\Delta_{w_0 \omega_i, \omega_i}(x) =  \Delta_{\omega_{i^*},w_0
\omega_{i^*}}(x) {\rm \ for \ all \ } i\} \ .$$ In particular, for
$G = SL_n$, we have
$$G^{w_0, w_0} = \{x \in G :
\Delta_{[n+1-i,n], [1,i]}(x) \neq 0, \Delta_{[1,i],[n+1-i,n]}(x)
\neq 0 {\rm \ for \ all \ } i\} \ ,$$ and
$$S^{w_0, w_0} = \{x \in G :
\Delta_{[n+1-i,n], [1,i]}(x) = \Delta_{[1,n-i],[i+1,n]}(x) \neq 0
{\rm \ for \ all \ } i\} \ ,$$ where $\Delta_{I,J} (x)$ is the
minor of a matrix $x$ with the row set $I$ and the column set $J$,
and $[a,b]$ stands for the set $\{a, a+1, \dots, b\}$. The
dimension of $S^{w_0, w_0}$ is equal to $2 \ell (w_0)$; for
$SL_n$, this amounts to $n(n-1)$.

To illustrate Corollary~\ref{cor:integrable-system}, consider the
symplectic leaf $S^{w_0, w_0}$ in $G = SL_3$. Its dimension is
$6$, and it is given in $G$ by the conditions
$$x_{31} =  \Delta_{[1,2],[2,3]}(x) \neq 0,\,\,  x_{13} = \Delta_{[2,3],[1,2]}(x) \neq 0\ ,$$
where $x_{ij}$ denotes the $(i,j)$-entry of the matrix $x$. Choose
the reduced word ${\bf j} = (1,2,1)$ of $w_0$. Using an explicit
description of the twist map $x \mapsto x'$ given in
\cite[Example~4.6]{fz}, the restrictions to $S^{w_0, w_0}$ of the
corresponding twisted minors $M_1, M_3$ and $M_5$ can be
calculated as follows:
\begin{align*}
&M_1 (x) = \Delta_{3,1} (x') = \frac{1}{x_{13}} \ ,\\
&M_3 (x) = \Delta_{[2,3],[1,2]}(x') = \frac{1}{x_{31}} \ ,\\
&M_5 (x) = \Delta_{2,2}(x') = \frac{x_{23}\Delta_{\{1,3\},
[1,2]}(x) - x_{13}\Delta_{[2,3],[1,2]}(x)}{x_{13}x_{31}} \ .
\end{align*}
By Corollary \ref{cor:integrable-system}, these functions form a
completely integrable system on $S^{w_0,w_0}$. In view of
Remark~\ref{remark}, the matrix entries $x_{31}, x_{13}$,
and $x_{22}$ also form a completely integrable system on $S^{w_0,w_0}$.}
\end{example}

\section{Proofs of the results on symplectic leaves}
\label{sec:proofs-symplectic}

We retain the terminology and notation of Section \ref{main-results}.
In particular, $G$ is a simply-connected connected semisimple complex Lie group
with the standard Poisson-Lie structure.
In this section, we prove Theorem \ref{symplectic-leaves},
Proposition~\ref{pr:Suv-toric-chart} and
Corollary~\ref{cor:number-of-components}.
Our starting point is the following description of
symplectic leaves in $G$.

\begin{proposition}
\label{pr:symplectic-leaves-thru-tildeS} The symplectic leaves in
$G$ are the connected components of the sets $\tilde S^{u,v} \cdot
a$ for some $u,v\in W$ and $a \in H$, where $\tilde S^{u,v}$ is
given by {\rm (\ref{eq:tildeS})}.
\end{proposition}

Proposition~\ref{pr:symplectic-leaves-thru-tildeS} appeared in
\cite{hl,r,y}. We still would like to outline a proof in order to
make the presentation more self-contained, and also since this
gives us a convenient occasion to introduce some notation needed
later.

\proof
We deduce Proposition~\ref{pr:symplectic-leaves-thru-tildeS} from
the following description of
symplectic leaves which is essentially due to M.~Semenov-Tyan-Shanski\u\i \ \cite{s-t}.
Let us identify $G$ with the diagonal of $G\times G$, and let
$G^*$ be the subgroup of $G\times G$ given by
$$
G^*=\{(b,b_-)\in B \times B_- : [b]_0=[b_-]_0^{-1}\}.
$$
Then the symplectic leaves of $G$ are the connected components of
the intersections of $G$ with the double cosets of $G^*$ in $G \times G$.

It remains to show that the intersections of $G$ with the double cosets of $G^*$ in $G \times G$
are precisely the sets $\tilde S^{u,v} \cdot a$
for $u,v\in W$ and $a \in H$.
By the definition, two elements $x$ and $y$ of $G$ belong to
the same double coset in $G^* \backslash (G \times G)/ G^*$
if and only if
\begin{equation}
\label{eq:double-coset}
y = b x b' = b_- x b'_-
\end{equation}
for some $b,b' \in B$ and $b_-,b'_- \in B_-$ such that
$[b]_0=[b_-]_0^{-1}$ and $[b']_0=[b'_-]_0^{-1}$.
In particular, $x$ and $y$ must belong to the same double Bruhat cell $G^{u,v}$.

To shorten the notation, let us define for $x\in G^{u,v}$
\begin{equation}
\label{eq:h-and-h'-1} h(x)=[{\overline{u}}^{-1}x]_0, \,\,
h'(x)=([x \overline{v^{-1}}]_0)^v.
\end{equation}
One easily checks that
\begin{equation}
\label{eq:h-and-h'-2} h(xa) = h(x) a, \,\, h'(xa) = h'(x) a, \,\,
h(ax)= a^uh(x), \,\, h'(ax)=a^vh'(x)
\end{equation}
for every $a\in H$.

Now suppose that $x, y \in G^{u,v}$ satisfy
(\ref{eq:double-coset}).
Using (\ref{eq:h-and-h'-2}), we obtain
$$h(y) =  h([b]_0 \cdot x \cdot [b']_0) =
[b]_0^u \cdot h(x) \cdot [b']_0,$$
and
$$h'(y) =  h'([b_-]_0 \cdot x \cdot  [b'_-]_0)
= [b_-]_0^v \cdot h'(x) \cdot [b'_-]_0 \ .$$
Therefore, $h(y) h'(y)$ and $h(x) h'(x)$ belong to the same coset
in $H/H^{u,v}$.

Conversely, suppose $h(y) h' (y)$ and $h(x) h'(x)$ belong to the
same coset in $H/H^{u,v}$, i.e. we have $a^v h(y) h'(y) = a^u h(x)
h'(x)$ for some $a \in H$. Setting
$$a' = a^v h'(y) h'(x)^{-1} = a^u h(x) h(y)^{-1},$$
we obtain
\begin{align*}
(x,x) &\in G^* (\overline u h(x), \overline{v^{-1}}^{\ -1} h'(x))
G^* = G^* (a \overline u h(x) (a')^{-1}, a^{-1}
\overline{v^{-1}}^{\ -1} h'(x)a') G^*\\
&= G^* (\overline u h(y), {\overline{v^{-1}}}^{\ -1} h'(y)) G^* =
G^* (y,y) G^* \ .
\end{align*}
Thus $x,y\in G^{u,v}$ satisfy
(\ref{eq:double-coset}) if and only if
\begin{equation}
\label{eq:same-coset}
\text{$h(x)h'(x)$ and $h(y)h'(y)$ belong to the same coset in $H/H^{u,v}$.}
\end{equation}
At the same time, (\ref{eq:tildeS}), (\ref{eq:h-and-h'-1}) and
(\ref{eq:h-and-h'-2}) imply that $x, y \in G^{u,v}$ belong to the
same set $\tilde S^{u,v} \cdot a$ if and only if
(\ref{eq:same-coset}) holds. This completes the proof of
Proposition~\ref{pr:symplectic-leaves-thru-tildeS}.
\endproof

Recall that $I(u,v)$ is the set of indices $i$ for which
$u\omega_i=v\omega_i=\omega_i$.
Let $\tilde H^{u,v}$ denote the subtorus of $H$ given by
$$
\tilde H^{u,v}= \{h\in H : h^{\omega_i}=1 \text{\ for \ } i\in I(u,v)\}.
$$
Clearly, $H^{u,v}$ is a subtorus of $\tilde H^{u,v}$.
In view of (\ref{eq:h-and-h'-1}), we have
\begin{align}
&\label{eq:tildeSuv} \tilde S^{u,v}=\{x\in G^{u,v}: h(x) h'(x) \in
H^{u,v}\} \ , \\
&\label{eq:Suv} S^{u,v}=\{x\in \tilde S^{u,v}: h(x)\in \tilde
H^{u,v}\} \ .
\end{align}

\begin{lemma}
\label{lem:h-over-h'}
For every $x \in G^{u,v}$, we have
$h(x) h'(x)^{-1} \in \tilde H^{u,v}$.
\end{lemma}

\proof
Let $i \in I(u,v)$.
Then we have
$$h(x)^{\omega_i} = ([{\overline{u}}^{\ -1}x]_0)^{\omega_i}
= \Delta_{u \omega_i, \omega_i} (x) = \Delta_{\omega_i, \omega_i}
(x) \ ,$$
and
$$h'(x)^{\omega_i} = ([x \overline{v^{-1}}]_0)^{v \omega_i}
= \Delta_{\omega_i, v^{-1} \omega_i} (x) =
\Delta_{\omega_i, \omega_i} (x) = h(x)^{\omega_i} \ ,$$
implying our statement.
\endproof

Let us denote
$$\Sigma=\{a \in H : a^2=e\} =
\{a \in H : a^{\omega_i} = \pm 1  \text{\ for \ all \ } i\} \ ,$$
and $\Sigma^{u,v}=\Sigma \cap \tilde H^{u,v}$. Thus, $\Sigma$ is a
finite Abelian group isomorphic to $(\ZZ/2\ZZ)^r$, and
$\Sigma^{u,v}$ is a subgroup of index $2^{|I(u,v)|}$ in $\Sigma$.
By (\ref{eq:h-and-h'-2}) and (\ref{eq:tildeSuv}), $\Sigma$ acts on $\tilde S^{u,v}$ by
right translations.
Using (\ref{eq:Suv}), we see that $\Sigma^{u,v}$ takes $S^{u,v}$ into
itself; furthermore, if $a$ and $a'$ belong to different cosets in
$\Sigma/\Sigma^{u,v}$ then $S^{u,v} \cdot a \cap S^{u,v} \cdot a'
= \emptyset$.
On the other hand, Lemma~\ref{lem:h-over-h'} and
(\ref{eq:tildeSuv}) imply that $h(x)^2 \in \tilde H^{u,v}$ for every $x\in
\tilde S^{u,v}$.
It follows that every $\Sigma$-orbit in $\tilde
S^{u,v}$ has non-empty intersection with $S^{u,v}$.
We conclude
that $\tilde S^{u,v}$ is the disjoint union of $2^{|I(u,v)|}$
right translates of $S^{u,v}$ by a set of representatives of
$\Sigma/ \Sigma^{u,v}$.
Since $S^{u,v}$ is closed in $\tilde
S^{u,v}$, it is also open, hence is a union of connected
components of $\tilde S^{u,v}$.

To complete the proofs of Theorem~\ref{symplectic-leaves}
and Corollary~\ref{cor:number-of-components}, it
remains to show that $S^{u,v}$ is connected.
This in turn is a consequence of Proposition~\ref{pr:Suv-toric-chart},
since the connectedness property of a complex algebraic variety
is preserved by passing to a dense Zariski open subset.
So it remains to prove Proposition~\ref{pr:Suv-toric-chart}.

In what follows, we fix $u,v \in W$ and a double reduced word
$\ii = (i_1, \dots, i_m)$ of $(u,v)$.
Recall that we append $r$ entries
$i_{m+1},\dots,i_{m+r}$ to $\ii$ by setting $i_{m+j}=\bar j$.
For any $j$, let $k(j)$ denote the
smallest index $k$ with $|i_k|=j$;
thus, $1 \leq k(j) \leq m$ for $j \notin I(u,v)$, and
$k(j) = m+j$ for $j \in I(u,v)$.

\begin{lemma}
\label{Suv-identification}
The variety $S^{u,v}$ can be identified
with the subvariety of pairs $(x,a)\in G^{u,v}\times H^{u,v}$
satisfying
\begin{equation}
\label{eq:definition-ofSuv-1}
M_{m+j}(x)=
\begin{cases}
1 & \text{if $j\in I(u,v) \,$;}\\
a^{-\omega_j}\prod_{i\notin
I(u,v)} M_{k(i)}(x)^{- \langle \alpha_i^\vee, v\omega_j \rangle} &
\text{if $j\notin I(u,v) \,$.}
\end{cases}
\end{equation}
\end{lemma}

\proof
By (\ref{eq:tildeSuv}) and (\ref{eq:Suv}),
$S^{u,v}$ can be identified with the subvariety of
pairs $(x,a)\in G^{u,v}\times H^{u,v}$ such that
$$h(x)\in \tilde H^{u,v},  \,\,
h(x)h'(x)=a \ ,$$
or, equivalently,
\begin{equation}
\label{eq:Suv-x-a}
h(x)^{\omega_j} = 1 \,\, (j \in I(u,v)), \quad
h(x)^{\omega_j} h'(x)^{\omega_j} = a^{\omega_j}
\,\, (j = 1, \dots, r) \ .
\end{equation}
It remains to show the equivalence of (\ref{eq:Suv-x-a}) and
(\ref{eq:definition-ofSuv-1}).

We claim that
\begin{equation}
\label{eq:h(x)-M}
h(x)^{\omega_j}= M_{m+j}(x)^{-1}, \,\,
h'(x)^{v^{-1}\omega_j}=
M_{k(j)}(x)^{-1} \
\end{equation}
for any $x \in G^{u,v}$ and $j = 1, \dots, r$.
Indeed, the definition (\ref{eq:twist}) of the twist
isomorphism $x \mapsto x'$ between $G^{u,v}$ and
$G^{u^{-1},v^{-1}}$ implies at once that
\begin{equation}
\label{eq:unbounded-minors-twist} [x' \overline v]_0 =
([{\overline u}^{\ -1} x]_0)^{-1}, \,\, [{\overline {u^{-1}}}^{\
-1} x']_0 = ([x \overline {v^{-1}}]_0)^{-1} \ .
\end{equation}
The equality (\ref{eq:h(x)-M})  is then a direct consequence of definitions
(\ref{eq:h-and-h'-1}) and (\ref{eq:twisted-minors-def}).

Since every fundamental weight $\omega_j$ can be written as
$\sum_i \langle \alpha_i^\vee, v\omega_j \rangle \cdot v^{-1} \omega_i$,
the equality (\ref{eq:h(x)-M}) implies that
$$h'(x)^{\omega_j} =
\prod_i M_{k(i)}(x)^{-\langle \alpha_i^\vee, v\omega_j \rangle} \ .$$
Substituting this and the first equality in (\ref{eq:h(x)-M}) into
(\ref{eq:Suv-x-a}), we conclude (after a slight simplification) that it is indeed
equivalent to (\ref{eq:definition-ofSuv-1}).
\endproof

Lemma \ref{Suv-identification} together with Theorem
\ref{biregular} immediately prove that the map
$x\mapsto(M_1(x),\dots, M_m(x), h(x)h'(x))$ provides an
isomorphism between $S^{u,v}\cap U_\ii$ and $(\mathbb
C-\{0\})^m\times H^{u,v}$.
To prove Proposition~\ref{pr:Suv-toric-chart}, it remains to show that
$S^{u,v}\cap U_\ii$ is dense in $S^{u,v}$.
Since $S^{u,v}$ is smooth, it
suffices to show that the complement of this set is of
(complex) codimension at least one.
In other words, we claim that $S^{u,v}\cap \{M_k=0\}$ is of codimension at
least one in $S^{u,v}$ for any $k = 1, \dots, m+r$.
But this is clear from the defining equations
(\ref{eq:definition-ofSuv-1}): the twisted minors $M_{m+j}$ and
$M_{k(i)}$ vanish nowhere on $S^{u,v}$, while the remaining $M_k$
just do not appear in (\ref{eq:definition-ofSuv-1}).
Proposition~\ref{pr:Suv-toric-chart},
Theorem \ref{symplectic-leaves} and Corollary
\ref{cor:number-of-components} are proved.

\section{Proofs of Theorem~\ref{poisson-brackets} and
Corollary \ref{cor:integrable-system}}
\label{sec:proofs-integrable}

In this section, we again fix two elements $u, v \in W$ and a double
reduced word $\ii = (i_1, \dots, i_m)$ of $(u,v)$. By
Theorem~\ref{biregular}, a generic element $x \in G^{u,v}$ has a
unique factorization $x = x_{\bf i}(a;t_1,\dots,t_m)$, so each
$t_k$ as well as $a^\gamma$ for any weight $\gamma$, can be viewed
as a rational function on $G^{u,v}$. The Poisson brackets between
these functions are given as follows.

\begin{lemma}
\label{t's}
We have
\begin{align}
&\{t_k,t_{k'}\}=\varepsilon(i_k)(\alpha_{|i_k|},\alpha_{|i_{k'}|})t_kt_{k'},
\text{ for }k<k'\\
&\{t_k,a^{\gamma}\}=(\alpha_{|i_k|},\gamma) t_k a^{\gamma}, \text{
for any weight }\gamma\\
&\{a^\gamma,a^{\gamma'}\}=0, \text{ for any weights
}\gamma,\gamma'.
\end{align}
\end{lemma}

A slight modification of this lemma can be found in~\cite{r+}. To
make the exposition more self-contained, we outline the proof.

\proof
Consider the Poisson-Lie group $SL_2^{(d)}$ (see
Example~\ref{example}), and two symplectic leaves in it given by
$$
S_+ = \{\mat{p}{q}{0}{p^{-1}}  :  p, q \neq 0\},  \quad
S_- = \{\mat{p}{0}{q}{p^{-1}}  :  p, q \neq 0\} \ ;
$$
in the notation of Theorem~\ref{symplectic-leaves}, we have
$S_+ = S^{e, w_0}$ and $S_- = S^{w_0, e}$.
In both these leaves, the Poisson structure is given by
$\{q,p\}= dpq$;
to indicate the dependence on $d$, we will write
$S_\pm = S_\pm^{(d)}$.

By the definition of the standard Poisson structure on $G$, a
double reduced word $\ii$ of $(u,v)$ gives rise to a Poisson map
$$\varphi_\ii: S_{{\rm sign}(\varepsilon (i_1))}^{(d_{|i_1|})} \times \cdots \times
S_{{\rm sign}(\varepsilon (i_m))}^{(d_{|i_m|})} \to G^{u,v}$$
given by
$$\varphi_\ii (g_1, \dots, g_m) = \varphi_{|i_1|} (g_1) \cdots
\varphi_{|i_m|} (g_m) \ .$$ We use the standard coordinates $(p_k,
q_k)$ in each factor $S_{{\rm sign}(\varepsilon
(i_k))}^{(d_{|i_k|})}$. An easy calculation using commutation
relations \cite[(2.5)]{fz} shows that
$$\varphi_{\ii}^*(a^\gamma)=\prod_{\ell=1}^m
p_\ell^{\langle \alpha^\vee_{|i_\ell|}, \gamma \rangle}, \,\,
\varphi_{\ii}^*(t_k)=q_k p_k^{-\varepsilon(i_k)} \prod_{\ell=k+1}^m p_\ell^{-\varepsilon(i_k) a_{|i_\ell|,|i_k|}} \ .$$
On the other hand, the only nonzero Poisson brackets between the
coordinates $p_1, q_1, \dots, p_m, q_m$ are $\{q_k, p_k\}=
d_{|i_k|}p_k q_k$. Recalling that $(\alpha_i, \gamma)=d_i \langle
\alpha_i^\vee, \gamma \rangle$, we obtain our statement.
\endproof

We will say that two functions $f,g$ on a
Poisson manifold \emph{quasi-commute} if $\{f,g\}=cfg$ for some
constant $c$.
In this situation, we denote $c=\langle f,g\rangle$.

Lemma~\ref{t's} asserts in particular that all functions $t_k$ and
$a^\gamma$ on $G^{u,v}$ quasi-commute with each other.
Clearly, any two monomials in quasi-commuting variables
quasi-commute, and the pairing $\langle f, g \rangle$ is
skew-symmetric and bilinear in the following sense:
$\langle f_1 f_2, g\rangle = \langle f_1 , g\rangle  +
\langle f_2, g\rangle$.
Taking into account (\ref{eq:M-k-monomial}),
we conclude that all twisted minors $M_k$ quasi-commute with each other.
To prove Theorem~\ref{poisson-brackets}, we only need to show that
\begin{equation}
\label{eq:MkMk'-pairing}
\langle M_k,M_{k'}\rangle= (\gamma^k,\gamma^{k'})-
(\delta^k,\delta^{k'})\text{ for }k<k'.
\end{equation}

We set
$$\gamma^k_\ell=u^{-1}_{\geq \ell} \gamma^k, \,\,
\delta^k_\ell=v^{-1}_{<\ell+1} \delta^k \ ;$$
in particular, we have
\begin{align*}
\gamma^k_1 = u \gamma^k, \,\, \gamma^k_k = \omega_{|i_k|}, \,\,
\delta^k_{k-1} = \omega_{|i_k|}, \,\, \delta^k_m = v^{-1} \delta^k\ .
\end{align*}
An easy inspection shows that we can
rewrite (\ref{eq:M-k-monomial}) as
\begin{equation}
\label{eq:Mk-cd}
M_k =a^{-u\gamma^k} \prod_{\ell=1}^{k-1} t_\ell^{c^k_\ell}
\prod_{\ell=k}^{m} t_\ell^{d^k_\ell} \ ,
\end{equation}
where the exponents $c^k_\ell$ and $d^k_\ell$ are determined from
\begin{equation}
\label{trail-property}
c^k_\ell\alpha_{|i_\ell|} = \gamma^k_{\ell}-\gamma^k_{\ell+1},
\,\, d^k_\ell\alpha_{|i_\ell|} = \delta^k_{\ell}-\delta^k_{\ell-1}
\ .
\end{equation}

To simplify calculations, we perform the following monomial change
of variables: replace each $t_k$ by
\begin{equation}
\label{u_k} y_k =
\begin{cases}
t_k & \text{if $\varepsilon(i_k)=+1 \,$;}\\
a^{-\alpha_{|i_k|}}t_k & \text{if $\varepsilon(i_k)=-1\,$.}
\end{cases}
\end{equation}
As an immediate consequence of Lemma~\ref{t's}, we get
\begin{align}
&\langle y_k,y_{k'}\rangle = \varepsilon(i_k)(\alpha_{|i_k|},\alpha_{|i_{k'}|})
\text{ for }k<k'\text{ and }\varepsilon(i_k)=\varepsilon(i_{k'}),\\
&\label{u_k,u_k'}\langle y_k,y_{k'}\rangle = 0 \text{ for }\varepsilon(i_k)\neq \varepsilon(i_{k'}),\\
&\langle y_k,a^{\gamma}\rangle = (\alpha_{|i_k|},\gamma) \text{ for any weight }\gamma,\\
&\langle a^\gamma,a^{\gamma'}\rangle =0 \text{ for any weights }\gamma,\gamma'.
\end{align}

Rewriting $M_k$ in these new variables, we get
$$M_k = a^{-\omega_{|i_k|}} \prod_{\ell=1}^{k-1}
y_\ell^{c^k_\ell} \prod_{\ell=k}^{m} y_\ell^{d^k_\ell} \ ;$$
the only thing to check here is that the weight
$(-u\gamma^k)$ in (\ref{eq:Mk-cd}) transforms into
$$-u \gamma^k + \sum_{\ell=1}^{k-1} c_\ell^k \alpha_{|i_\ell|} =
-u \gamma^k + \sum_{\ell=1}^{k-1}
(\gamma^k_{\ell}-\gamma^k_{\ell+1}) = -u \gamma^k +
\gamma^k_1-\gamma^k_k = -\omega_{|i_k|} \ .$$

Now everything is ready for the proof of (\ref{eq:MkMk'-pairing}).
We have
\begin{align*}
\langle M_k,M_{k'}\rangle &= \langle a^{-\omega_{|i_k|}},
M_{k'}\rangle +\langle M_k,a^{-\omega_{|i_{k'}|}}\rangle
+\langle\prod^{k-1}_{\ell=1}y_\ell^{c^k_\ell},\prod^{k'-1}_{\ell'=1}y_{\ell'}^{c^{k'}_{\ell'}}\rangle\\
&+\langle\prod^m_{\ell=k}y_\ell^{d^k_\ell},\prod^m_{\ell'=k'}y_{\ell'}^{d^{k'}_{\ell'}}\rangle
+\langle\prod^{k-1}_{\ell=1}y_\ell^{c^k_\ell},\prod^{m}_{\ell'=k'}y_{\ell'}^{d^{k'}_{\ell'}}\rangle
+\langle\prod^m_{\ell=k}y_\ell^{d^k_\ell},\prod^{k'-1}_{\ell'=1}y_{\ell'}^{c^{k'}_{\ell'}}\rangle\\
\end{align*}
The last two terms vanish by (\ref{u_k,u_k'}). Let us calculate
the remaining terms. First of all, we have
\begin{align*}
\langle a^{-\omega_{|i_k|}}, M_{k'}\rangle&=\langle
a^{-\omega_{|i_k|}},\prod^{k'-1}_{\ell=1}y_{\ell}^{c^{k'}_{\ell}}\prod^{m}_{\ell=k'}
y_{\ell}^{d^{k'}_{\ell}} \rangle
=(\omega_{|i_k|},\sum^{k'-1}_{\ell=1}c^{k'}_{\ell}\alpha_{|i_{\ell}|}
+\sum^{m}_{\ell=k'}d^{k'}_{\ell}\alpha_{|i_{\ell}|})\\
&=(\omega_{|i_k|},\sum^{k'-1}_{\ell=1}
(\gamma^{k'}_{\ell}-\gamma^{k'}_{\ell+1})
+\sum^{m}_{\ell=k'}(\delta^{k'}_{\ell}-\delta^{k'}_{\ell-1}))\\
&=(\omega_{|i_k|},\gamma^{k'}_1-\gamma^{k'}_{k'}+
\delta^{k'}_m-\delta^{k'}_{k'-1} )
=(\omega_{|i_k|},u \gamma^{k'}+ v^{-1}
\delta^{k'}-2\omega_{|i_{k'}|}) \ .
\end{align*}
Similarly,
\begin{equation*}
\langle M_k,a^{-\omega_{|i_{k'}|}}\rangle=
-(u \gamma^{k}+ v^{-1} \delta^{k}-2\omega_{|i_{k}|},\omega_{|i_{k'}|}) \ .
\end{equation*}
The third term can be calculated as follows:
\begin{align*}
\langle\prod^{k-1}_{\ell=1}y_\ell^{c^k_\ell},\prod^{k'-1}_{\ell'=1}y_{\ell'}^{c^{k'}_{\ell'}}\rangle
&=\sum^{k-1}_{\ell=1}(c^k_\ell\alpha_{|i_\ell|}, \sum^{\ell-1}_{\ell'=1}c^{k'}_{\ell'}\alpha_{|i_{\ell'}|}
- \sum^{k'-1}_{\ell'=\ell+1}c^{k'}_{\ell'}\alpha_{|i_{\ell'}|})\\
&=\sum^{k-1}_{\ell=1}(\gamma^k_\ell-\gamma^k_{\ell+1},\sum^{\ell-1}_{\ell'=1}(\gamma^{k'}_{\ell'}-\gamma^{k'}_{\ell'+1})
-\sum^{k'-1}_{\ell'=\ell+1}(\gamma^{k'}_{\ell'}-\gamma^{k'}_{\ell'+1}))\\
&=\sum^{k-1}_{\ell=1}(\gamma^k_\ell-\gamma^k_{\ell+1},
\gamma^{k'}_1  + \gamma^{k'}_{k'} - \gamma^{k'}_\ell - \gamma^{k'}_{\ell+1})\\
&=\sum^{k-1}_{\ell=1}(\gamma^k_\ell-\gamma^k_{\ell+1},u \gamma^{k'}+\omega_{|i_{k'}|})
-\sum^{k-1}_{\ell=1}(\gamma^k_\ell-\gamma^k_{\ell+1},\gamma^{k'}_{\ell}+\gamma^{k'}_{\ell+1})\\
&=(\gamma^k_1-\gamma^k_k,u \gamma^{k'}+\omega_{|i_{k'}|})-
\sum^{k-1}_{\ell=1}((\gamma^k_\ell,\gamma^{k'}_{\ell})-
(\gamma^k_{\ell+1},\gamma^{k'}_{\ell+1}))\\
&=(u \gamma^k-\omega_{|i_k|},u \gamma^{k'}+\omega_{|i_{k'}|})-
(u \gamma^k,u \gamma^{k'})+(\omega_{|i_k|},\gamma^{k'}_k) \\
&= - (\omega_{|i_k|},u \gamma^{k'}) + (u \gamma^k, \omega_{|i_{k'}|})
- (\omega_{|i_k|}, \omega_{|i_{k'}|}) + (\omega_{|i_k|},\gamma^{k'}_k) \ .
\end{align*}
Here we used the fact that, for $\ell<k$ and
$\varepsilon(i_\ell)=-1$, we have
$(\gamma^k_\ell,\gamma^{k'}_{\ell+1})=
(s_{|i_\ell|}\gamma^k_{\ell+1},\gamma^{k'}_{\ell+1})=
(\gamma^k_{\ell+1},\gamma^{k'}_{\ell})$.

The last remaining term is calculated in the same way.
We leave the details to the reader and only give the answer:
$$\langle\prod^m_{\ell=k}y_\ell^{d^k_\ell},\prod^m_{\ell'=k'}y_{\ell'}^{d^{k'}_{\ell'}}\rangle
= - (\omega_{|i_k|}, v^{-1} \delta^{k'}) + (v^{-1} \delta^k, \omega_{|i_{k'}|})
+ (\omega_{|i_k|}, \omega_{|i_{k'}|}) - (\delta^{k}_{k'-1}, \omega_{|i_{k'}|}) \ .
$$

Combining all the terms together and performing numerous
cancellations, we obtain:
$$\langle M_k,M_{k'}\rangle = (\omega_{|i_k|},\gamma^{k'}_k)
-(\delta^k_{k'-1},\omega_{|i_{k'}|}) \ .$$
Since $(\omega_{|i_k|},\gamma^{k'}_k)=(\gamma^k,\gamma^{k'})$ and
$(\delta^k_{k'-1},\omega_{|i_{k'}|})=(\delta^k,\delta^{k'})$, this
completes the proof of Theorem~\ref{poisson-brackets}. \hfill $\Box$

\smallskip

\noindent{\it Proof of Corollary~\ref{cor:integrable-system}.}
Since the torus $H^{u,u}$ is trivial, Proposition~\ref{pr:Suv-toric-chart}
implies that the twisted minors $M_1,\dots, M_m$
form a system of local coordinates on $S^{u,u}$;
in particular, the functions $M_{2k-1}$ are algebraically independent on~$S^{u,u}$.
Moreover, there are exactly $\frac{1}{2}\dim S^{u,u}$ of them.

It remains to show that the functions $M_{2k-1}$ pairwise Poisson
commute on~$S^{u,u}$. This follows at once from
Theorem~\ref{poisson-brackets} and the following computation for
$k<k'$:
\begin{align*}
&(\gamma^{2k-1},\gamma^{2k'-1})=
(s_{j_{\ell(u)}}\cdots s_{j_k} \omega_{j_k},
s_{j_{\ell(u)}}\cdots s_{j_{k'}} \omega_{j_{k'}})\\
&=(s_{j_{k'-1}}\cdots s_{j_k}\omega_{j_k},\omega_{j_{k'}})
=(\omega_{j_k},  s_{j_k} \cdots s_{j_{k'-1}} \omega_{j_{k'}})\\
&=(s_{j_1}\cdots s_{j_{k-1}} \omega_{j_k},
s_{j_1}\cdots s_{j_{k'-1}} \omega_{j_{k'}})
= (\delta^{2k-1},\delta^{2k'-1}).
\end{align*}
Corollary~\ref{cor:integrable-system} is proved.  \hfill $\Box$

\end{document}